\documentclass[12pt, oneside]{article}

\usepackage[dvips]{graphicx}
\usepackage{amsthm}
\usepackage{amsmath}
\usepackage{amsfonts}
\usepackage{amssymb}

\theoremstyle{plain}
\newtheorem{thm}{\textbf{Theorem}}[section]
\newtheorem{lem}[thm]{Lemma}
\newtheorem{defn}[thm]{Definition}

\newcommand{\ignore}[1]{}
\begin{document}
\title{Cover Pebbling Numbers and Bounds for Certain Families of Graphs}
\author{Nathaniel G. Watson \\Washington University in St. Louis \\ \ttfamily ngwatson@wustl.edu \\ \\ Carl R. Yerger \\ Harvey Mudd College \\ \ttfamily cyerger@hmc.edu}
\maketitle

\abstract{Given a configuration of pebbles on the vertices of a graph, a
\emph{pebbling move} is defined by removing two pebbles from some vertex
and placing one pebble on an adjacent vertex.  The cover pebbling number of
a graph, $\gamma(G)$, is the smallest number of pebbles such that through a
sequence of pebbling moves, a pebble can eventually be placed on every
vertex simultaneously, no matter how the pebbles are initially distributed.
The cover pebbling number for complete multipartite graphs and wheel graphs
is determined. We also prove a sharp bound for $\gamma(G)$ given the
diameter and number of vertices of $G$. }\newline \newline \noindent

\section{Introduction}
\phantom{space }One recent development in graph theory, suggested by
Lagarias and Saks, called pebbling, has been the subject of much research
and substantive generalizations.  It was first introduced into the
literature by Chung\cite{Chung}, and has been developed by many others
including Hurlbert, who published a survey of pebbling results in
\cite{Glenn2}. Given a connected graph $G$, distribute $k$ pebbles on its
vertices in some configuration, $C$.  Specifically, a \emph{configuration} on a
graph $G$ is a function from $V(G)$ to $\mathbb{N} \cup \{0\} $
representing an arrangement of pebbles on $G.$ We call the total number of pebbles $k$ the \emph{size} of the configuration.
A \emph{pebbling move} is defined as the simultaneous removal of two pebbles from some vertex and addition of one pebble on an
adjacent vertex.  A pebble can be moved to a root vertex $v$ if it is
possible to place one pebble on $v$ in a sequence of pebbling moves. We
define the pebbling number, $\pi(G)$ to be the minimum number of pebbles
needed so that for any initial distribution of pebbles, it is possible to
move to any root vertex $v$ in $G$.

The concept of cover solvability was introduced in \cite{Glenn1}. We call a
configuration on a graph \emph{cover solvable} if, starting with this
configuration, it is possible, through a sequence of pebbling moves, to
simultaneously place one pebble on every vertex of the graph. The
\emph{cover pebbling number} of a graph, $\gamma(G)$, is defined as the
smallest number such that every configuration of this size is cover
solvable. One application in \cite{Glenn1} for $\gamma(G)$ is based on a
military application where troops must be distributed simultaneously.

In \cite{Glenn1}, the cover pebbling number for complete graphs,
paths and trees is determined. Also, Hurlbert and Munyan
\cite{Glenn3} are preparing a proof that determines the cover pebbling
number of the $d$-cube.

This paper will consider various questions related to cover pebbling,
including open problem 9 in \cite{Glenn1}.   Section 2 describes the
computation for calculating the cover pebbling number for complete
multipartite graphs.  In Section 3, we will compute the cover pebbling
number of $W_n$, the wheel graph.   We conclude the paper in Section 4 by
constructing a tight upper bound for the cover pebbling number of graphs
with specified diameter $d$ and number of vertices $n$.

\section{Complete Multipartite Graphs}
\begin{defn}
For $ s_1 \geq s_2 \geq \cdots \geq s_r$, let $K_{s_1, s_2, \ldots, s_r}$
be the complete $r$-partite graph with $s_1,s_2,\ldots,s_r$ vertices in
vertex classes $c_1, c_2, \ldots, c_r$ respectively.
\end{defn}

\begin{defn}
For a complete $r$-partite graph $G=K_{s_1, s_2, \ldots, s_r},$ let
$\phi(G)=4s_1 + 2s_2 + \cdots + 2s_r - 3.$
\end{defn}

\begin{thm}
$\gamma(K_{s_1, s_2, \ldots, s_r}) = \phi(G). $
\end{thm}

\begin{proof}
First, we will show that not every configuration of size $\phi(K_{s_1, s_2,
\cdots, s_r}) - 1 $ on $K_{s_1, s_2, \cdots, s_r}$ is cover-solvable.
Consider the case where all $\phi(K_{s_1, s_2, \cdots, s_r}) - 1$ pebbles
are on one vertex of $c_1$, call it $x$. There are $k=s_2+s_3+\cdots+s_r$
vertices that are distance $1$ from $x$ and $l=s_1-1$ vertices that are
distance $2$ from $x$.  For the $k$ vertices a distance $1$ from $x$, $2k$
pebbles are required to cover these vertices, and for the $l$ vertices at
distance $2$ from $x$, there are $4l$ pebbles required to cover these
vertices. We need one more pebble to remain on $x$, for a total of
$2k+4l+1=\phi(K_{s_1, s_2, \cdots, s_r})$ pebbles required, which is one
more than we have. Thus, this configuration is not cover-solvable.

Now suppose that there exists some complete $r$-partite graph $K_{s_1, s_2,
\ldots, s_r}$ which has a configuration of size $\phi(K_{s_1, s_2,
\ldots, s_r})$ that is not cover-solvable. Among such graphs, choose one of minimal order
(let it be $G'= K_{s_1',s_2',\ldots, s_{r'}'}).$

First, we will show that $G'$ cannot be a star graph (that is, a
$K_{s_1',1}$.) To see this, consider:

\begin{defn}
\cite{Glenn1} Let $T$ be a tree and let $V(T)$ be the vertex set of $T$.  For $v \in V(T)$, define
$$s(v) = \sum_{u \in V(T)} 2^{d(u,v)},$$
where $d(u,v)$ denotes the distance from $u$ to $v$, and let
$$s(T) = \max_{v \in V(T)} s(v).$$
\end{defn}
In \cite{Glenn1} it is shown that for any tree $T$, $\gamma(T) = s(T)$.
Since $G'$ is a tree, we can compute $\gamma(G')$ by evaluating $s(v)$ for
all $v \in G$ to obtain $s(G)$.  If $v \in c_1$ then $s(v) = 4s_1'-1$, and
if $v \in c_2$ then $s(v) = 2s_1' + 1$.  Thus, $s(G') = \gamma(G') =
4s_1'-1 = 4s_1' + 2s_2' - 3=\phi(G')$.  Hence, for a star, every
configuration of size $\phi(G')$ is cover-solvable.  Since $G'$ is not a
star, further suppose that for any $G'$, each complete multipartite
subgraph $G$ of $G'$ is cover-solvable with $\phi(G)$ pebbles.

Notice that for any complete $p$-partite graph with $p \geq 2$ other than a
star graph, the removal of a vertex from the graph leaves a subgraph that
is a complete $q$-partite graph with $q \geq 2.$  Since $G'$ cannot be a
star, for any vertex $v \in G',$ $G'-v$ is a complete $r^*$-partite graph
with $r^* \geq 2.$ Furthermore, since by our assumption of the minimality of $G'$, for any complete
$r$-partite graph $G$ smaller than $G',$ a configuration of size
$\phi(G)$ or greater must be cover-solvable, and since clearly
$\phi(G'-v) \leq \phi(G')-2,$ any configuration of size $\phi(G')-2$
or greater on $G'-v$ is cover-solvable.

Let $C$ be a configuration of size $\phi(G')$ on $G'$. Suppose $C(v)=1$
or $2$ for some $v \in G'.$ Then $C$ restricted to $G'-v$ is a
configuration of size at least $\phi(G')-2$ and thus is cover-solvable on
$G'-v.$ After we carry out the steps of the cover-solution of this
subgraph, we will have cover-solved $G'$, contradicting our hypothesis.

Otherwise, if $C(v)=0$ or $C(v)\geq 3$ for all $v \in G'$, choose some $v'$
for which $C(v')=0$ (if no such $v'$ exists, we are done).  Then consider
the vertices of $G'$ which are in different vertex classes of $G'$ from
$v'.$ If at least one of these is initially occupied, call it $v''.$ Then
since $C(v'') \geq 3$, we can cover $v'$ with pebbles from $v''$, while
leaving $\phi(G')-2$ pebbles on $G'-v'.$ Thus, the configuration of
pebbles on $G'$ after this move, restricted to the subgraph $G'-v'$ is
cover solvable, and after we carry out the steps of the cover-solution of
this subgraph, we will have cover-solved $G'$. Otherwise, all the vertices
in the vertex classes of $G'$ that are different than the vertex class from
$v'$ are empty. Thus, all pebbles are on vertices in the vertex class of
$v',$ and in particular, some vertex $w$ of this class has pebbles on it,
so $C(w) \geq 2$. Thus, we can use pebbles on $w$ to cover some vertex $w'$
in another vertex class, as all these vertices are empty. Note that after
this move, the configuration of pebbles on $G'-w'$ has size $\phi(G')-2$,
and thus this configuration restricted to the subgraph $G'-w'$ is
cover-solvable. Again, after we carry out the steps of the cover-solution
of this subgraph, we will have cover-solved $G'$.

\end{proof}

\section{The Wheel Graph}
In this section, we will compute $\gamma(W_n)$, where $W_n$ is the wheel
graph.  The wheel graph is composed of a cycle consisting of $n$ vertices,
$v_1, \ldots, v_n$, which are all connected to a hub vertex, $v_0$, for a
total of $v = n+1$ vertices.
\begin{thm}
For $n \geq 3$, $\gamma(W_n) = 4n - 5 = 4v - 9$.
\end{thm}

\begin{proof}
Consider the configuration of pebbles where all the pebbles are on one
vertex of $W_n$, say $x$, that is not the hub.  In this case, $2$ pebbles
are required to cover each of the three vertices adjacent to $x$, and $4$
pebbles are required to cover each of the $n-3$ vertices that are a
distance of $2$ away from $x$.  The total number of pebbles required to
cover-solve these vertices is $4n-6$. However, we require one more pebble
to place on $x$.  Hence, $\gamma(W_n) \geq 4n-5$.

To complete the proof, we will show that if there is some configuration of
pebbles on $W_n$ with at least $4n-5$ pebbles, then the configuration is
cover-solvable.  Suppose $C$ is a configuration of pebbles on $W_n$ and
consists of at least $4n-5$ pebbles.  We now will describe a sequence of
moves that will cover-solve any such configuration.  First, if there are
outer vertices on $W_n$ that are empty but are adjacent to outer vertices
$w$ such that initially $C(w) \geq 3$, then if the adjacent vertices can be
covered and $w$ can also remained covered, then these adjacent vertices
should be covered. Let $k$ be the number of outer vertices that are covered
after this process.

\textbf{Case 1:} Suppose that $k = 0$.  In this case, all the pebbles are
on the hub vertex.  To cover-solve the remaining $v-1$ vertices, we can
cover $\lfloor\frac{4v-10}{2}\rfloor = 2v-5$ vertices using the excess
pebbles already on the hub vertex.  Since $v \geq 4$ and $2v-5 \geq v-1$,
we can cover-solve all of the outer vertices in this manner.

\textbf{Case 2:}  Suppose that $k = 1$ or $k = 2$.  Each outer vertex
covered in the process above requires at most two pebbles to cover it.
Since $v \geq 4$, there are at least $4v-9-2k$ pebbles already on the hub
vertex. After subtracting $1$ pebble for the hub itself, there are
$4v-10-2k$ pebbles that can be used such that pebbles can be placed on the
remaining $v-k-1$ uncovered vertices.  With these remaining pebbles on the
hub, we can cover at least $\lfloor \frac{4v-10-2k}{2} \rfloor = 2v-5-k$
vertices. Since $2v-5-k \geq v-k-1$ for $v\geq 4$, there are enough pebbles
to cover-solve $W_n$ in this situation.

\textbf{Case 3:}  Suppose that $k \geq 3$.  Again, each outer vertex in the
process above requires at most two pebbles to cover it.  If there are any
pairs of pebbles remaining on outer vertices such that removing the pairs
would not uncover that vertex, those pairs of pebbles should be moved to
the hub vertex.  After this process, there are at least $\lceil \frac{4v -
9 - 2k}{2} \rceil = 2v - 4 - k$ pebbles on the hub vertex.  Notice that
this bound is based on the worst case that occurs when no pebbles are
initially on the hub vertex.  From the hub vertex, it takes exactly $2$
pebbles to cover each of the remaining outer vertices and one pebble to
cover the hub vertex. So at most $\lfloor \frac{2v - 5 - k}{2} \rfloor = v
- 3 - \lfloor\frac{k}{2}\rfloor$ outer vertices can be pebbled. Since there
are at most $v-k-1$ outer vertices left to be pebbled, and for $k \geq 3$,
$v-k-1 \geq v - 3 - \lfloor\frac{k}{2}\rfloor$, there are enough pebbles to
cover-solve $W_n$ in this case, and the proof is complete.

\end{proof}

\section{The Cover Pebbling Number of Graphs of Diameter $d$}
\begin{defn}
A \emph{binary weighting} on a graph $G$ is a function from $V(G)$ to
$\{0,1\}.$ If $B$ is a binary weighting on $G,$ then let the order $|B|$ of
$B$ be $\sum_{v \in G} B(v).$
\end{defn}
\begin{defn}
For a graph $G$ and binary weighting $B$, a configuration $C$ on $G$ will
be called \emph{permissible} (with respect to $B$) if for all $v \in G,$
$B(v)=0 \implies C(v)=0.$ A permissible configuration on a graph $G$ with a
binary weighting $B$ will be called \emph{cover-solvable} (with respect to
$B$) if we can reach a configuration on which $B(v)=1 \implies C(v) \geq 1$
for all $v \in G$ by a sequence of pebbling moves.
\end{defn}

\begin{lem} Let $G$ be a graph of diameter $d$, $B$ a binary weighting on $G,$ and $C$
a configuration of size at least $(|B|-1)2^d + 1$ on $G$ which is
permissible with respect to $B.$ Then $C$ is cover-solvable with respect to
$B.$
\end{lem}
\begin{proof}
Assume the opposite. Then for all pairs $\{G,B\}$ of a graph $G$ together
with a binary weighting on $G$ such that there exists a non-cover-solvable
configuration of size at least $(|B|-1)2^d + 1$ (where $d$ is the diameter
of $G$,) choose one for which $|B|$ is minimal, and call it $\{G',B'\}.$
Let $d'$ be the diameter of $G'$, let $k =(|B'|-1)2^{d'} + 1,$ and choose
some configuration (call it $C'$) on $G'$ which is permissible with respect
to $B'$, has size at least $k$ and is not cover-solvable.

Certainly we cannot have $|B'|=1$, for then the only permissible
configuration of size $|C'| \geq k = 1$ is the function which takes the
value $|C'|$ on the lone vertex for which $B'=1$, and is zero elsewhere.
This configuration covers all vertices with non-zero weights, and so is
trivially cover-solvable, creating a contradiction.

Now, suppose that $|B'| \geq 2.$ If it is true that $C'(v) > 0 $ whenever
$B'(v) = 1$ we have a contradiction, for $C'$ is then trivially
cover-solvable. Otherwise, let $v'$ be some vertex of $G'$ for which
$C'(v') = 0$ and $B'(v')=1.$ At most $|B'|-1$ vertices are initially
occupied, and there are at least  $=(|B'|-1)2^{d'} + 1$ total pebbles, so
by the pigeonhole principle, there are at least $2^{d'}+1$ pebbles on some
vertex (call it $v''$). Since the diameter of $G'$ is $d',$ $d(v',v'') \leq
d'.$ Thus we can move $2^{d'}$ of the pebbles from $v''$ onto $v',$ through
a series of pebbling moves, losing half of these pebbles for each edge we
must move across, but leaving at least one pebble on $v'$ if we move all
pebbles via one of the shortest paths.

Now, define a binary weighting $B^*$ on $G$ by $$B^*(v)=\left\{
\begin{array}{r@{\quad:\quad}l} 0 & v =v'
\\ B'(v) & v \not= v'\end{array} \right.$$
and define a configuration $C^*$ on $G$ by $$C^*(v)=\left\{
\begin{array}{r@{\quad:\quad}l} 0 & v = v' \\ C'(v'')-2^{d'} & v=v'' \\
C'(v) & \textnormal{otherwise} \end{array} \right.$$ This is the
configuration after we have moved pebbles from $v''$ onto $v',$ except that
we ignore the pebbles on $v'$ and designate it as a vertex which need not
be covered by pebbles. Clearly $|B^*| = |B'|-1$ and $|C^*|=|C'|-2^{d'}$ so
from $|C'| \geq ((|B'|-1)2^{d'} + 1),$ we see $|C^*| \geq ((|B^*|-1)2^{d'}
+ 1).$ $C^*$ is permissible with respect to $B^*,$ and so by our assumption
of the minimality of $B'$, $C^*$ is cover-solvable with respect to $B^*.$

If we carry out the moves of the cover-solution of $C^*$ on $G$ starting
with the configuration left on $G'$ after our initial movement of pebbles
from $v''$ to $v'$, (certainly this is possible because this configuration
is no smaller than $C^*$ on any vertex,) we will have covered every vertex
of $G'$ for which $B^*=1.$ Also, we must still have $v' \geq 1$, because
$C^*(v')=0$, which does not permit any sequence of moves that decreases the
number of pebbles on $v'.$ Thus every vertex for which $B'=1$ now has $C'
\geq 1,$ and we have cover-solved $C'$ with respect to $B',$ which
contradicts the assumption that $C'$ was not cover-solvable.

\end{proof}

\begin{thm} Let $G$ be a graph of order $n$ and diameter $d$, and let $C$ be a configuration
on $G$ of size at least $2^d(n-d+1)-1.$ Then $G$ is cover-solvable
\textnormal{(}with respect to the weighting on $G$ which is equal to $1$
for each vertex.\textnormal{)}
\end{thm}
\begin{proof}

First, notice that this bound is sharp because we can exhibit the following
class of graphs where $\gamma(G) \geq 2^d(n-d+1)-1.$  Suppose we have a
graph consisting of $n$ vertices and diameter $d$.  Then we will construct
a fuse graph (a path connected to a star) whose length is $d-1$ and has
$n-d$ spokes at the end of the fuse. Here is an example for $n = 7$ and $ d
= 4$.

\begin{figure}[htb]
\unitlength 1mm
\begin{center}
\begin{picture}(60,25)
\put(05,15){\circle*{3}} \put(15,15){\circle*{3}} \put(25,15){\circle*{3}}
\put(35,15){\circle*{3}} \put(45,15){\circle*{3}} \put(45,05){\circle*{3}}
\put(45,25){\circle*{3}}
\put(05,15){\line(1,0){10}} 
\put(15,15){\line(1,0){10}} \put(25,15){\line(1,0){10}}
\put(35,15){\line(1,0){10}} \put(35,15){\line(1,1){10}}
\put(35,15){\line(1,-1){10}}
\end{picture}
\end{center}
\caption{A graph where $n = 7$ and $d=4$ such that $\gamma(G) =$
$2^d(n-d+1)-1.$} \label{ex2}
\end{figure}
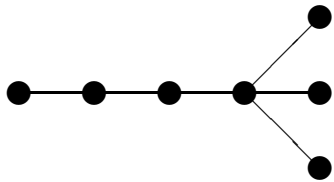

Suppose all the pebbles are on the last vertex of the path, which is at
distance $d$ from the spokes (in the figure, the leftmost vertex.) Then
each of the $n-d$ spokes requires $2^d$ pebbles, and the path requires $2^d
- 1$ vertices to cover-solve. (Note: In $\cite{Glenn1},$ the cover-pebbling
number of all trees is found. Thus, we know for these particular trees that
$\gamma(G) = 2^d(n-d+1)-1$ even before proving this theorem.)

We will prove this theorem by defining an algorithm by induction which will
take us to a configuration, the solvability of which we can prove using the
lemma. Let $R_0 = \{ v \in G : C(v)>0\},$ let $S_0 = \{ v \in G :
C(v)=0\},$ and let $T_0 = \emptyset.$ Let $C_0  = C.$

For illustrative purposes, we now describe the first step of the algorithm.
If $S_0=\emptyset,$ we are clearly done, for $C$ already covers $G$.
Otherwise, note that since $R_0$ and $S_0$ are complementary, there exist
vertices $r_0 \in R_0$ and $s_0 \in S_0$ such that $d(r_0,s_0) = 1.$ If
$C_0(r_0)=1$ or $C_0(r_0)=2,$ then let $R_1 = R_0\setminus\{r_0\},$
$S_1=S_0$ and  $T_1 = T_0\cup\{r_0\}=\{r_0\}.$ In this case, let $C_1 =
C_0.$

If on the other hand $C_0(r_0) \geq 3$ then we move $2$ pebbles from $r_0$
to $s_0,$ and instead put $s_0$ in $T_1$ and define $C_1$ according to the
following configuration. Explicitly, in this case let $R_1=R_0$, $S_1 = S_0
\setminus\{ s_0\},$ and $T_1 = T_0\cup\{s_0\}=\{s_0\}.$ Define $C_1$ on $G$
by $$C_1(v) =\left\{
\begin{array}{r@{\quad:\quad}l} r_0-2 & v = r_0 \\ 1 & v=s_0 \\ C_0(v) &
\textnormal{otherwise} \end{array} \right.$$

Define the sequences $R_0, R_1, \ldots, R_{d-1},$ $S_0, S_1, \ldots,
S_{d-1},$ $T_0, T_1, \ldots, T_{d-1},$ and $C_0, C_1, \ldots, C_{d-1},$
recursively in an analogous manner. Suppose for some $m < d-1$ we have
$R_m,$ $S_m,$ $T_m$, and $C_m,$ such that the following hold:
\begin{enumerate}
\item $|T_m| = m.$
\item $R_m,$ $S_m$ and $T_m$ are disjoint and $R_m \cup S_m \cup T_m = V(G).$
\item For all $v \in R_m \cup T_n, \ \ C_m(v) > 0$ and for all $v \in S_m, \ \ C_m(v)=0.$
\item $C_m$ is a configuration which can be reached from $C$ by a sequence of pebbling moves.
\item $R_m$ and $S_m$ are both non-empty.
\item $\sum_{v \in R_m} C_m(v) \geq [2^d(n-d+1)-1]-[2^{m+1}-2].$
\end{enumerate}
Note that all these conditions are trivially true for $m=0$.

From condition 1, it is evident that the minimum distance between $R_m$ and
$S_m$ is at most $m+1.$ Take points $r_m \in R_m$ and $s_m \in S_m$ for
which this minimum distance is achieved (and thus $d(r_m , s_m) \leq m+1.
$) If $C_m(r_m) \leq 2^{m+1}$ then let $R_{m+1} = R_m \setminus \{r_m\},$
$S_{m+1}=S_m$ and  $T_{m+1} = T_m\cup\{r_m\}.$ In this case, let $C_{m+1} =
C_m.$

Otherwise, if $C_m(r_m) > 2^{m+1}$ then we can move $2^{m+1}$ pebbles along
a minimal path from $r_m$ to $s_m,$ which is of length at most $m+1.$ We
lose half of these pebbles for each edge we must move across, but we will
be able to move $2^{(m+1)-d(r_m, s_m)} \geq 1$ onto $s_m.$ Put $s_m$ in
$T_{m+1}$ and define $C_{m+1}$ according to the configuration after these
moves. Explicitly, in this case let $R_{m+1}=R_m,$ $S_{m+1} = S_m
\setminus\{ s_m\}$ and $T_{m+1} = T_m\cup\{s_m\}.$ Define $C_{m+1}$ on $G$
by $$C_{m+1}(v) =\left\{
\begin{array}{r@{\quad:\quad}l} r_m-2^{m+1} & v = r_m \\ 2^{(m+1)-d(r_m,
s_m)} & v=s_m \\ C_m(v) & \textnormal{otherwise} \end{array} \right.$$

For $m+1,$ it is clear from our definitions that conditions 1, 2, 3, and 4
still hold. Condition 6 also holds, for in either of the two above cases,
the total number of pebbles left on $R_{m+1}$ is at most $2^{m+1}$ less
than were on $R_m$. Thus, \begin{eqnarray*} \sum_{v \in R_{m+1}} C_{m+1}(v)
 & \geq &  \sum_{v \in R_m} C_m(v)-2^{m+1} \\ & \geq &
[2^d(n-d+1)-1]-[2^{m+1}-2]-2^{m+1}  \nonumber \\ & = &
[2^d(n-d+1)-1]-[2^{m+2}-2]. \nonumber \end{eqnarray*}  For condition 5,
since always $m+1 < d$ and $n \geq d,$ $[2^d(n-d+1)-1]-[2^{m+1}-2] > 0.$
Thus, the fact that condition 6 is true for $m+1$ necessitates that $R_m
\not= \emptyset.$ Also, if $S_{m+1} = \emptyset$ then $C_{m+1}(v)>0$ for
all $v \in R_m \cup S_m \cup T_m = V(G),$  and since $C_m$ is attainable
from $C$ by a sequence of pebbling moves, we have cover-solved $C$ and we
are done. So we may assume $S_{m+1} \not= \emptyset$ and condition 5 holds.

By this recursive definition, we now have $R_{d-1},$ $S_{d-1},$ $T_{d-1},$
and $C_{d-1}$ for which conditions 1-6 hold. Now define a binary weighting
$B$ on $G$ by
$$B(v)=\left\{ \begin{array}{r@{\quad:\quad}l} 1 & v \in R_{d-1} \cup S_{d-1} \\ 0 & v \in T_{d-1}\end{array} \right.$$
Also, define $C_{d-1}'$ on $G$ by
$$C_{d-1}'(v)=\left\{ \begin{array}{r@{\quad:\quad}l} C_{d-1}(v) & v \in R_{d-1} \cup S_{d-1} \\ 0 & v \in T_{d-1}\end{array} \right.$$ Clearly
$C_{d-1}'$ is permissible with respect to $B.$  From condition 1,
$|T_{d-1}|=d-1$ so $|B| = n-d+1$, and from condition 6 we have $|C_{d-1}'|
\geq [2^d(n-d+1)-1]-[2^{(d-1)+1}-2]=2^d(n-d)+1.$  Thus, by the lemma,
$C_{d-1}'$ is cover-solvable with respect to $B$.

By condition 4, $C_{d-1}$ is a configuration which can be reached from $C$
by a sequence of pebbling moves. If after we carry out this sequence of
moves, we carry out the moves of this cover-solution of $C_{d-1}'$ on $G$
(certainly this is possible because $C_{d-1}'$ is no greater than $C_{d-1}$
on any vertex,) we will have covered every vertex of $G$ for which $B=1,$
that is every vertex in $R_{d-1} \cup S_{d-1}.$ Also, every vertex $v \in
T_{d-1}$ must remain covered, because for each of these vertices,
$C_{d-1}'(v)=0,$ which does not permit any sequence of moves which
decreases the number of pebbles on $v.$ Applying, condition 2, we see for
every vertex $v \in V(G)=R_{d-1} \cup S_{d-1} \cup T_{d-1},$ our final
configuration after this sequence of moves is greater than zero, and so we
have cover-solved $C$.

\end{proof}

\paragraph*{Acknowledgment}

The authors received support from NSF grant DMS-0139286, and would like to
acknowledge East Tennessee State University REU director Anant Godbole for
his guidance and encouragement.

\end{document}